\documentclass[11pt]{amsart}
\begin{document}

\numberwithin{equation}{section}

\def\R{{\mathbb R}}
\def\K{{\mathcal K}} 
\def\L{{\mathcal L}}
\def\O{{\Omega}}
\def\Tr{\mbox{\rm{Tr}}}

\def\vep{{\varepsilon}}
\def\p{{\partial}}
\def\a{{\alpha}}
\def\b{{\beta}}
\def\e{{\eta}}
\def\g{{\gamma}}
\def\s{{\sigma}}
\def\t{{\tau}}
\def\E{{\mathbb E}}
\def\N{{\mathbb N}}
\def\F{{\mathcal F}}
\def\P{{\mathbb P}}
\def\EE{{\mathcal E}}
\def\o{{\omega}}
\def\l{{\ell}}
\def\vp{{\varphi}}
\def\un{{\mathrm{1~\hspace{-1.4ex}l}}}
\def\to{{\rightarrow}}
\def\th{\theta}
\def\Dt{{\Delta t}}
\def\ds{\displaystyle}
\def\be{\begin{equation}}
\def\ee{\end{equation}}
\def\ba{\begin{array}}
\def\ea{\end{array}}
\def\l{\left}
\def\r{\right}
\def\D{{{\mathbb D}^{1,2}}}
\newcommand{\Frac}{\displaystyle \frac}
\newcommand{\Sum}{\displaystyle \sum}
\newcommand{\Int}{\displaystyle \int}
\newcommand{\Sup}{\displaystyle \sup}

\newtheorem{Theorem}{Theorem}[section]
\newtheorem{Definition}[Theorem]{Definition}
\newtheorem{Proposition}[Theorem]{Proposition}
\newtheorem{Lemma}[Theorem]{Lemma}
\newtheorem{Corollary}[Theorem]{Corollary}
\newtheorem{Remark}[Theorem]{Remark}
\newtheorem{Example}[Theorem]{Example}
\newtheorem{Hypothesis}[Theorem]{Hypothesis}

\title[Weak approximation of SPDEs]
{Weak approximation of stochastic partial differential equations: the non linear case}

\author[A. Debussche]{
Arnaud Debussche}
\thanks{ENS de Cachan, Antenne de Bretagne,
Campus de Ker Lann, Av. R. Schuman,
35170 BRUZ, FRANCE ({\tt arnaud.debussche@bretagne.ens-cachan.fr})}

\thanks{Acknowledgments: Part of this work was done while the author visited the 
Institut Mittag-Leffler (Djursholm, Sweden) 
during the semester "Stochastic Partial Differential Equations".}

\keywords{Weak order, stochastic heat equation, Euler scheme.}

\subjclass{
35A40, 60H15, 60H35.
}

\begin{abstract}

We study the error of  the Euler scheme applied to a stochastic partial differential 
equation. We prove that as it is often the case, the weak order of convergence is twice
the strong order. A key ingredient in our proof is Malliavin calculus which enables 
us to get rid of the irregular terms of the error. We apply our method to the 
case a semilinear stochastic heat equation driven by a space-time white noise. 

\end{abstract}

\maketitle

\today

\section{Introduction}

When one 
considers a numerical scheme for a stochastic equation,  two types of errors can be considered.
The strong error measures the pathwise approximation of the true solution by a numerical one.
This problem has been extensively studied in finite dimension for stochastic differential 
equations (see for instance \cite{kloeden-platten},  \cite{milstein}, \cite{milstein-tretyakov}, \cite{talay}) and also more recently
in infinite dimension for various types of stochastic partial differential equations (SPDEs) (see  
among others \cite{allen-novosel-zhang}, \cite{dBDschema}, \cite{davie-gaine}, 
\cite{greksch-kloeden}, \cite{gyongy}, \cite{gyongy2},
\cite{gyongy-millet}, \cite{gyongy-millet2}, \cite{gyongy-millet3}, 
\cite{gyongy-nualart}, \cite{hausenblas}, \cite{hausenblas2}, \cite{lord-rougemont}, \cite{millet-morien}, 
\cite{Printems}, \cite{shardlow}, \cite{yan1}, \cite{yan2}, \cite{walsh}). Another way to measure the error
is the so-called  weak order of convergence of a numerical scheme which is concerned with
the approximation of the law of the solution at a fixed time. In many applications, this error is more relevant.
Pioneering work by Milstein (\cite{milstein-1}, \cite{milstein-2})
 and Talay (\cite{talay2}) have been followed by many articles (see references in the books cited above).
Very few works exist in the literature for the weak approximation of solution of SPDEs. A delayed stochastic 
differential equation has been studied in \cite{buckwar-shardlow}. 
Weak order for a SPDE has been studied only recently in \cite{dB-D}, \cite{DP}, \cite{haus3}. In order to 
explain the novelty of the present article, let us focus on a specific example.
 
We consider a stochastic nonlinear heat equation in a bounded interval $I=(a,b)\subset R$
with Dirichlet boundary conditions and driven by a space-time white noise:
\be
\label{e0.1}
\left\{
\ba{l}
\ds \frac{\partial X}{\partial t} = X_{\xi\xi} +f(X) +\s(X) \dot \eta,\; \xi\in I,\; t>0,\\
\\
\ds X(a,t)=X(b,t)=0,\; t>0,\\
\\
\ds X(\xi,0)=x(\xi),\,  \xi\in I.
\ea
\right.
\ee
Where $f$ and $\s$ are smooth Lipschitz functions from $\R$ to $\R$. 

We introduce the classical abstract framework extensively used in the book 
\cite{DPZ}. We set $H=L^2(I)$, $A=\partial_{\xi\xi}$, $D(A)=H^2(I)\cap H^1_0(I)$,
$W$ is a cylindrical Wiener process so that the space-time white noise is mathematically represented as
the time derivative of $W$. We set $f(x)(\xi)=f(x(\xi))$, $x\in H$ and define
$\s\; :\; H\to \L(H)$ by $\s(x)h (\xi) = \s(x(\xi))h(\xi)$, $x,h\in H$. We then rewrite \eqref{e0.1} as
\begin{equation}
\label{e0.2}
 \left\{
\begin{array}{l}
dX = (AX +f(X))dt + \sigma(X) dW,\\
X(0)=x.
\end{array}
\right.
\end{equation}
It is well known that this equation has a unique solution.
We investigate the error committed when approximating this solution by the solution
of the Euler scheme
\begin{equation}
\label{e0.3}
 \left\{
\begin{array}{l}
X_{k+1}-X_k= \Dt\l(AX_{k+1} +f(X_k)\r) + \s(X_k)\l( W((k+1)\Dt) -W(k\Dt)\r),\\
X_0=x,
\end{array}
\right.
\end{equation}
where $\Dt=T/N$, $N\in\N$, $T>0$. 

The study of the weak error aims to prove bounds of the type:
$$
\l|\E(\varphi(X(n\Dt))) -\E(\varphi(X_n)\r| \le c \Dt^\delta,
$$
with a constant $c$ which may depend on $\varphi,\, x,\, N$ and on the various parameter in the equation.
Also $\varphi$ is assumed to be a smooth function on $H$. If such a bound is true, we say that the scheme
has weak order $\delta$. In comparison, the strong error is given by 
$\E(| (X(n\Dt))- X_n|)$ or $\E(\sup_{n=0,\dots,N} | (X(n\Dt))- X_n|)$. Clearly, if the scheme has strong 
order $\tilde \delta$ then it has weak order $\delta\ge \tilde\delta$. Indeed, the test functions
$\varphi$ are Lipschitz. In general, it is expected that the weak order is larger than the the strong order.

In the case of the Euler scheme applied to a stochastic differential equation, it is well known that the
strong order is $1/2$ whereas the weak order is $1$ (see \cite{talay}). The classical proof of this uses 
the Kolmogorov equation associated to the stochastic equation. The main difficulty to generalize this proof 
to the infinite dimensional equation \eqref{e0.2} is that this Kolmogorov equation is then a partial differential
equation with an infinite number of variables and involving unbounded operators (see \eqref{e3.4} below).
The delayed stochastic differential equation studied in \cite{buckwar-shardlow} is an infinite 
dimensional problem but since the equation does not contain differential operators the Kolmogorov 
equation is simpler to study. In \cite{haus3}, a SPDE similar to \eqref{e0.2} is considered but very 
particular test functions $\varphi$  are used. They are allowed to depend only on finite dimensional projections of the unknown and the bound of the  weak error involves a constant which strongly depends on the 
dimension. In \cite{dB-D}, \cite{DP}, the Kolmogorov equation is not used directly. A change of variable 
is used in order to simplify it. In \cite{dB-D}, the stochastic nonlinear Schr\"odinger equation is considered
and the fact that the linear Schr\"odinger equation generates an invertible group is used in an essential 
way. This is obviously wrong for the heat equation considered here. The same change of unknown 
works in the case of a linear equation with additive noise as shown in \cite{DP} but there it is used that
the solution can be written down explicitly. We have not been able to generalize this idea to the 
non linear equation considered here. 

We use in fact the original method developed by Talay in the finite dimensional case. The weak error
is decomposed thanks to the Kolmogorov equations on each time step. Each term represents the error 
between the solution of the Kolmogorov equation on one time step and the approximation given by
the numerical solution. Due to the presence of unbounded operators, this apparently requires a lot
of smoothness on the numerical solution. The main idea here is to observe that the non smooth
part of the solutions of \eqref{e0.2} and \eqref{e0.3} are contained in a stochastic integral. We get 
rid of this stochastic integral thanks to Malliavin calculus and an integration by part.
We are thus able to prove that as expected the weak order is twice the strong order without 
artificial assumption except from a technical one on $\s$. We restrict our presentation to the abstract 
equation above, a nonlinear heat equation driven by a space-time white noise. However, our 
method is general and can be used for more general equations as will be shown in future articles. 
Also, we only consider a semi-discretization in time. A full discretization will treated in forthcoming works.

Note that the method developed here does allow to recover the result of \cite{DP}. Indeed, in
the Euler scheme \eqref{e0.3}, the linear term is fully implicit and we cannot consider a scheme
where it is partially implicit such as the theta-scheme considered in \cite{DP}. Note also that the
proof below are much more complicated than in \cite{DP} and \cite{dB-D}. 

Malliavin calculus has already been used for the numerical analysis of stochastic equations.
In \cite{bally-talay}, it is used to prove an expansion of the error of the Euler scheme for a stochastic differential equation under minimal assumptions on the test functions $\varphi$. This is a completely 
different idea and the Malliavin calculus is used completely differently. It is not clear that such 
ideas could be used for a SPDE.  In a different spirit, Malliavin calculus is used in \cite{Szepessy-tempone-zouraris}
to analyse adaptive schemes for the weak approximation of stochastic differential equations.

Our method is much closer to the method developped in \cite{KohatsuHiga}. There, the Malliavin 
calculus is also used to get rid of a stochastic integral which appears when writting 
down the weak error. However, it is done in a global way and the error is not decomposed 
as in the present  article. A fundamental feature of Kohatsu-Higa's method is that the Kolmogorov 
equation is
not used so that more general stochastic equation can be can considered. The solution does not
need to be markovian. However, no SPDE have been considered with this method.

\section{Preliminaries and main result}

We consider the following stochastic partial differential equation written in an abstract form 
in a Hilbert space $H$ with norm $|\cdot|$ and inner product $(\cdot,\cdot)$:
\begin{equation}
\label{e2.1}
 \left\{
\begin{array}{l}
dX = (AX +f(X))dt + \sigma(X) dW,\\
X(0)=x,
\end{array}
\right.
\end{equation}
where the unknown $X$ is a random process on a probability space $(\O,\F,\P)$ depending
on $t>0$ and on the initial data $x\in H$. The operator $A$ is a negative 
self-adjoint operator on $H$ with domain $D(A)$ and has a compact inverse. We assume that
\be
\label{e2.2}
\Tr ((-A)^{-\a}) <\infty, \mbox{ for all } \alpha>1/2.
\ee
We define classically the domain $D((-A)^{\beta })$, $\beta \in R$, of fractional powers of $A$ and
set
$$
|x|_\beta  = |(-A)^{\beta }x|, \; x\in D((-A)^\beta) .
$$ 
The nonlinear function $f$ takes values in $H$ and is assumed to be $C^3$ with bounded 
derivatives up to order $3$. We denote by $L_f$ a constant such that for $x,y\in H$
\be
\label{e2.3}
\ba{l}
|f(x)|\le L_f (|x|+1),\\
|f(x)-f(y)|\le L_f |x-y|,\\
|f'(x)-f'(y)|_{\L(H)}\le L_f |x-y|.
\ea
\ee
The noise is written in terms of a cylindrical Wiener process $W$ on $H$ (see \cite{DPZ}) associated 
to a filtration $(\F_t)_{t\ge 0}$. The nonlinear mapping acting on the noise maps $H$ onto
$\L(H)$, it  is also assumed to be $C^3$ with bounded 
derivatives up to order $3$. We denote by $L_\s$ a constant satisfying
\be
\label{e2.4}
\ba{l}
|\s(x)|_{\L(H)}\le L_\s (|x|+1),\\
|\s(x)-\s(y)|_{\L(H)}\le L_\s |x-y|.\\
\ea
\ee
We need a stronger assumption on this mapping, we require
\be
\label{e2.4bis}
|\s''(x)\cdot(h,h)|_{\L(H))}\le L_\s |h|^2_{-1/4},\; x\in H,\; h\in H.
\ee 
Note that this implies a strong restriction on $\s$. (See Remark \ref{r1} below for some 
comments on this assumptions).

Recall that the cylindrical Wiener process can be written as 
$$
W=\sum_{\ell\in\N} \beta_\ell e_\ell
$$
where here and in the following $(e_\ell)_{\ell\in\N}$ is any orthonormal basis of $H$ and 
$(\beta_\ell)_{\ell\in\N}$ is an associated sequence of independent brownian motions. 
This series does not converge in $H$ but in any larger Hilbert space $U$ such
that the embedding $H\subset U$ is Hilbert-Schmidt. Similarly, given a linear operator $\Phi$
from $H$ to a possibly different Hilbert space $K$, the Wiener process 
$\Phi W=\sum_{\ell\in\N} \beta_\ell \Phi e_\ell$ is well defined in $K$ provided 
$\Phi\in \L_2(H,K)$, the space of Hilbert-Schmidt operators from $H$ to $K$. (See the definition 
just below).

Recall also that the stochastic integral $\int_0^T \Psi(s)dW(s)$ is defined as an element 
of $K$ provided that $\Psi$ is an adapted process with values in $\L_2(H,K)$ such that
$\int_0^T |\psi(s)|_{\L_2(H,K)}^2 ds <\infty$ {\it a.s.} (see \cite{DPZ}).

If $L \in {\L}(H)$ is a nuclear operator, $\Tr(L)$ denotes the trace of the operator $L$, i.e. 
$$
\Tr(L) = \sum_{i\geq 1} (L e_i,e_i) < + \infty.
$$
It is well known that the previous definition does not depend on the choice of the Hilbertian 
basis. Moreover, the following properties hold for $L$ nuclear and $M$ bounded
\begin{equation} \label{eq:trace_prop1}
\Tr{}(LM) = \Tr{}(ML), 
\end{equation}
\noindent and, if $L$ is also positive,
\begin{equation} \label{eq:trace_prop2}
\Tr{}(LM) \leq \Tr{}(L) \|M\|_{{\L}(H)}.
\end{equation}

Hilbert-Schmidt operators play also an important role. An operator $L \in {\L}(H)$ is Hilbert-Schmidt if $L^*L$ is a  nuclear operator on $H$. We denote by ${\L}_2(H)$ the space of
such operators. It is a Hilbert space for the norm
$$
\|L\|_{{\L}_2(H)}= \left( \Tr{} (L^*L)\right)^{1/2}=  \left( \Tr{} (L L^*)\right)^{1/2}.
$$
It is classical that  if $L\in {\L}_2(H)$,
$M\in {\L}(H)$, $N\in {\L}(H)$ then $NLM\in {\L}_2(H)$ and 
\begin{equation} \label{eq:HS}
\| NLM\|_{{\L}_2(H)}\le \| N\|_{{\L}(H)} 
\|L\|_{{\L}_2(H)}\|M\|_{{\L}(H)}.
\end{equation}
See \cite{DPZ}, appendix C, or \cite{gohberg-krejn} for more details on nuclear and Hilbert-Schmidt 
operators. Note that \eqref{e2.2} implies that $(-A)^{-\beta}$ is Hilbert-Schmidt for any $\beta>1/4$.

\bigskip

Our assumptions imply that for any $x\in H$, there exists a unique solution $X(t)$ to equation to 
\eqref{e2.1} (see for instance \cite{DPZ}, chapter 7). 
In the sequel, we often recall the dependence of the solution on the initial data 
by using the notation $X(t,x)$.

We approximate equation \eqref{e2.1} by an implicit Euler schemes. Let $\Dt=\frac{T}N>0$ be a 
time step, we define the sequence $(X_k)_{k=0,\dots,N}$ by 
\be
\label{e2.5}
\l\{
\ba{l}
X_{k+1}= S_\Dt X_k +\Dt S_\Dt f(X_k) + \sqrt\Dt S_\Dt \s(X_k) \chi_{k+1},\\
X_0=x.
\ea
\r.
\ee
We have set $\chi_{k+1} = (W((k+1)\Dt)-W(k\Dt))/\sqrt\Dt$.
The operators $S_\Dt$ is defined by
$$
S_\Dt=(I-\Dt A)^{-1}.
$$
This is the classical fully implicit Euler scheme. 
It will be convenient to use the integral form of \eqref{e2.1}
\be
\label{e2.6}
X(t)=S(t)x+\int_0^t S(t-s) f(X(s))ds + \int_0^t S(t-s) \s(X(s))dW(s),\; t\ge 0,
\ee
where $S(t)=e^{tA}$ is the semigroup generated by $A$. 
Similarly, \eqref{e2.5} can be rewritten as 
\be
\label{e2.7}
X_k= S_\Dt^k x +\Dt \sum_{\ell=0}^{k-1} S_\Dt^{k-\ell}  f(X_\ell) + 
\sqrt\Dt\sum_{\ell=0}^{k-1} S_\Dt^{k-\ell}  \s(X_{\ell}) \chi_{\ell+1}.
\ee
It will be convenient in the following to use the notation:
$$
t_k=k\Dt,\; k=0,\dots,N.
$$
The following inequalities are classical and easily proved using the spectral decomposition of 
$A$

\be
\label{e1}
\l|(-A)^\beta  S_\Dt^k \r|_{\L(H)}\le c t_k^{-\beta },\quad k\ge 1,\quad \beta  \in [0,1].
\ee

\be
\label{e1bis}
\l|(-A)^\beta  S(t) \r|_{\L(H)}\le c t^{-\beta },\quad t>0,\quad \beta  \ge 0.
\ee

\be
\label{e3}
\l|(-A)^\beta   S_\Dt \r|_{\L(H)}\le c \Dt^{-\beta },\quad \beta  \in [0,1].
\ee

\be
\label{e4}
\l|(-A)^{-\beta }\l(  I-S_\Dt\r) \r|_{\L(H)}\le c \Dt^{\beta },\quad \beta  \in [0,1].
\ee

Note that in \eqref{e2.7} and \eqref{e2.5}, the noise term  makes sense in $H$. Indeed, by 
\eqref{e3}, \eqref{e2.2} and \eqref{eq:HS}, 
we know that $S_\Dt$ is a Hilbert-Schmidt operator on $H$.

We are interested in the approximation of the law of the solution of \eqref{e2.1}. More precisely,
we wish to prove an estimate on the error committed when approximating $\E(\varphi(X(T,x)))$ by
$\E(\varphi(X_N(x)))$. The function $\varphi$ is a smooth function on $H$. 

In all the article, we use the notation $D\varphi(x)$ for the differential of a $C^1$ function on $H$
at the point $x$. If $\varphi:\; H\mapsto K$, where $K$ is another Hilbert space, $D\varphi(x)
\in \L(H,K)$ the space of continuous linear operator from $H$ to $K$. When $K=\R$, we identify 
the differential with the gradient thanks to Riesz identification theorem. We use the same notation and 
have the identity for $x,h\in H$: 
$$
D\varphi(x).h=(D\varphi(x),h).
$$
Similarly, if $\varphi\in C^2(H,\R)$, $D^2\varphi(x)$ is a bilinear operator from $H\times H$ 
to $\R$ and can be identified with a linear operator on $H$ through the identity:
$$
D^2\varphi(x).(h,k)=(D^2\varphi(x)h,k),\; x,h,k\in H.
$$
Sometimes, we also use the notations $\varphi'$, $\varphi''$ instead of $D\varphi$ or 
$D^2\varphi$.

Given two Banach spaces $K_1$ and $K_2$, we denote by $\|\cdot\|_k$ the norm on 
$C_b^k(K_1,K_2)$, the space of $k$ times continuously differentiable mapping from 
$K_1$ to $K_2$ with derivatives bounded up to order $k$.

We use Malliavin calculus in the course of the proof. We now recall the basic definitions.
(See \cite{nualart}). 
Given a smooth real valued function $F$ on $H^n$ and $\psi_1,\dots,\psi_n\in L^2(0,T,H)$, 
the Malliavin derivative of
the smooth random variable 
$F(\int_0^T (\psi_1(s),dW(s)),\dots,\int_0^T (\psi_n(s),dW(s)))$ at time $s$ in the direction 
$h\in H$ is given by
$$
\ba{l}
\ds D_s^h \l[ F\l(\int_0^T (\psi_1(s),dW(s)),\dots,\int_0^T (\psi_n(s),dW(s))\r) \r]\\
\\
\ds =\ds\sum_{i=1}^n \partial_i F\l(\int_0^T (\psi_1(s),dW(s)),\dots,\int_0^T (\psi_n(s),dW(s))\r) 
\l(\psi_i(s),h\r). 
\ea
$$
We also define the process $DF$ by $(DF(s),h)=D^h_sF$.
It can be shown that $D$ defines a closable operator with values in $L^2(\O,L^2(0,T,H))$ and 
we denote by $\D$ the closure of the set of smooth random variables as above for the topology 
defined by the norm 
$$
\|F\|_\D= \l(\E(|F|^2)+\E(\int_0^T |D_sF|^2ds\r)^{1/2}.
$$
We define similarly the Malliavin derivative of random variables taking values in $H$. 
If $G=\sum_{i\in\N} F_i e_i \in L^2(\O,H)$ where $F_i\in \D$ for all $i\in\N$ and 
$\sum_{i\in\N} \int_0^T |D_sF_i|^2ds <\infty$, we set 
$D^h_s G=\sum_{i\in \N} D^h_sF_i e_i$, $D_s G=\sum_{i\in \N} D_sF_i e_i$. 
We define $\D(H)$ as the set of such random variables.

When $h=e_m$, we write $D^{e_m}=D^m$.

The chain rule is valid and given $u\in C_b^1(\R)$, $F\in\D$ then $u(F)\in \D $
and $D(u(F))=u'(F)DF$. Also if $G=\sum_{i\in\N} F_i e_i\in\D(H)$ and $u\in C^1_b(H,\R)$ then
$u(G)\in\D$ and 
$D(u(G))= Du(G).DG=(Du,DG)$, or equivalently 
$D_s^h(u(G))= \sum_{i\in\N}\partial_i u D^h_s F_i=(Du,D^h_sG)$. 

Note that as already mentionned, we identify the differential of a function in $C^1(H,\R)$ with its 
gradient.

For $F\in\D$ and $\psi\in L^2(\O\times[0,T];H)$ such 
that $\psi(t)\in \D$ for all $t\in [0,T]$ and $\int_0^T\int_0^T |D_s \psi(t)|^2 ds dt <\infty$, we 
have the integration by part formula:
$$
\E\l(F\int_0^T (\psi(s),dW(s))\r)= \E\l(\int_0^T (D_sF,\psi(s))ds\r),
$$
where the stochastic integral is a Skohorod integral which is in fact defined by duality.
In this article, we only need to consider the Skohorod integral of adapted processes in which 
case it corresponds with the It\^o integral. Moreover, the integration by part formula above
holds for $F\in\D$ and $\psi\in L^2(\O\times[0,T];H)$ when $\psi$ is an adapted process. Recall that if $F$ is $\F_t$ measurable then 
$D_sF=0$ for $s\ge t$.

We will often use the following form of the integration by part formula whose proof is left to the reader.
\begin{Lemma}
\label{l2.1}
Let $F\in\D(H)$, $u\in C^2_b(H)$ and $\psi\in L^2(\O\times [0,T],\L_2(H))$ be an adapted 
process
then 
$$
\ba{ll}
\ds \E\l( Du(F)\cdot\int_0^T\psi(s)dW(s) \r)&\ds =\E\l( \sum_{m\in\N}\int_0^T D^2u(F)\cdot( D_s^mF,\psi(s)e_m)ds\r)\\
\\
&\ds = \E\l( \int_0^T \Tr \l( \psi^*(s) D^2u(F)D_sF \r) ds \r).
\ea
$$
\end{Lemma}

Also we remark that this Lemma remains valid if $u$ is not assumed to be bounded 
but only  $u\in C^2(H)$ provided the expectations and the integral above are well
defined. This is easily seen by approximation of $u$ by bounded functions.

We now state our main result.

\begin{Theorem}
\label{t1.1}
Assume that $f$ and $\s$ are $C^3_b$ functions from $H$ to $H$ and $\L(H)$ and that $\s$ satisfies 
\eqref{e2.4}, then for any $x\in H$, $T>0$, $\vep>0$, the Euler Scheme \eqref{e2.5} satisfies the 
following weak error estimate
$$
|\E(\varphi(X(T,x))) - \E(\varphi(X_N))| \le C(T,|\varphi|_{C^3_b}, |x|, \vep) \Dt^{1/2-\vep} ,
\quad \varphi\in C^3_b(H).
$$
\end{Theorem}

\begin{Remark}
\label{r1}
Assumption \eqref{e2.4} is quite restrictive. It is void for an additive noise or a noise of the
form $BX\,dW$ where $B$ is a linear operator from $H$ to $\L(H)$. Otherwise, it implies that
the noise is a perturbation of such noise. An example of a noise satisfying this is
$$
\s(x)=Bx + \tilde\s((-A)^{-1/4}x)
$$
where $B\in \L(H)$ and $\tilde \s \,:\, H \to \L(H)$ is a $C^3$ function with derivatives bounded up to order $3$.
This assumption is crucial in our proof. It is used in essential 
way in Lemma \ref{l1bis} which is used at many points of the proof.  

Apart from this point, our result is 
optimal. If the noise is assumed to satisfied some non degeneracy assumptions, the smothness 
assumption on 
the test function $\varphi$ can be weakened. This will be investigated in a future work.
\end{Remark}

In all the article, $C$ or $c$ denote  constants which may depend on $A,f,\s, Q$ or $T$ but not
on $\Dt$. Their 
value may change from one line to another. The initial data $x$ is fixed and the constant
may also depend on $|x|$. Note also that we assume that $\Dt \le 1$, we could also assume 
$\Dt\le \Dt_0$ for some $\Dt_0>0$. In this case, the different constants would depend on 
$\Dt_0$. Finally, $\vep$ is a small positive number.

\section{Proof of the main result}

The proof uses different tools from stochastic calculus such as It\^o formula, 
Kolmogorov equations, Malliavin calculus. Sometimes, it may be very 
lengthy and technical to justify rigorously their use in infinite dimension. 
We avoid these tedious justifications by using Gakerkin approximations. 
We replace 
equation \eqref{e2.1} by the finite dimensional stochastic equation
$$
dX_m = (AX_m + f_m(X_m)) dt + \s_m(X_m) dW, \quad X_m(0)=P_m
$$
where $P_m$ is the eigenprojector on the $m$ first eignevectors of $A$, 
$f_m(x) =P_m f(x)$, $\s_m(x)=P_m\s(x)P_m$. It is not difficult to prove that
$X_m$ converges to $X$ in various senses. 

Similarly, we replace the discrete unknown $X_k$ by a finite dimensional 
sequence defined in an obvious way.

We prove the result for these finite dimensional objects with constants that
do not depend on the dimension $m$. It is then easy to deduce the result for our 
infinite dimensional equation.

In order to lighten the notation, we omit to explicit the dependence on $m$ below
and write
$X$, $f$, $\s$ instead of $X_m$, $f_m$, $\s_m$.  

\bigskip 

{\bf Step 1:} We first define a continuous interpolation of the discrete unknown. 

\smallskip

We rewrite \eqref{e2.5} as follows:
$$
X_{k+1}= X_k +  \int_{t_k}^{t_{k+1}} A_\Dt X_k+  S_\Dt f(X_k) ds +
\int_{t_k}^{t_{k+1}}   S_\Dt \s(X_k)  dW(s) 
$$
where $A_\Dt= S_\Dt A$. Note that 
$A_\Dt$ is in fact a Yosida regularization of $A$ and is a bounded operator: 
\be
\label{e3.0}
|A_\Dt|_{\L(H)} \le c \Dt^{-1}.
\ee
It is then natural to define $\tilde X$ on $[0,T]$ by
\be
\label{e3.0bis}
\tilde X (t) =  X_k +  \int_{t_k}^{t} A_\Dt X_k+  S_\Dt f(X_k) ds +
  \int_{t_k}^{t}   S_\Dt \s(X_k)  dW(s),\quad t\in [t_k,t_{k+1}). 
\ee
Clearly, $\tilde X$ is a continuous and adapted process. Given a smooth 
function $G$ on $[0,T]\times H$, 
It\^o formula implies for $t\in [t_k,t_{k+1})$ (see \cite{DPZ}):
\be
\label{e3.1}
\ba{ll}
G(t,\tilde X(t))& \ds = G(t_k,\tilde X(t_k)) + \int_{t_k}^t \frac{dG}{dt} (s,\tilde X(s)) +
 L_{k,\Dt} G (s,\tilde X(s)) ds \\
\\
&\ds+ \int_{t_k}^t (DG(s,\tilde X(s)), \s(X_k)dW(s)).
\ea
\ee
Where for $\psi \in C^2(H,\R)$
$$
L_{k,\Dt} \psi (x)= \frac12 \Tr \l\{  \l( S_\Dt \s(X_k) \r)\l( S_\Dt \s(X_k) \r)^*  D^2 \psi(x) \r\}
+(A_\Dt X_k + S_\Dt f(X_k), D\psi(x)).
$$

\bigskip

{\bf Step 2:} Decomposition of the error.

\smallskip

Let us define 
\be
\label{e3.2}
u(t,x)=\E(\varphi (X(t,x))),\; t\in [0,T].
\ee
Then the weak error at time $T$ is equal to 
\be
\label{e3.3}
\ba{ll}
u(T,x)-\E(\varphi(X_N)) &= \E(u(T,x) -u(0,X_N)\\
& =\ds \sum_{k=0}^{N-1} \E\l(u(T-t_k,X_k)-u(T-t_{k+1},X_{k+1})\r).
\ea
\ee
It is well known that $u$ is a solution to the forward Kolmogorov equation:
\be
\label{e3.4}
\ba{ll}
\ds \frac{du}{dt}(t,x)&= Lu(t,x)\\
\\
&=\ds  \frac12 \Tr\{ \s(x) \s^*(x) D^2u(t,x) \} +(Ax+f(x),Du(t,x)).
\ea
\ee
Therefore, It\^o formula \eqref{e3.1} implies
$$
\E(u(T-t_{k+1}, X_{k+1})) = \E(u(T-t_k,X_k)) + \E\int_{t_k}^{t_{k+1}} L_{k,\Dt} u(T-t,\tilde X(t)) - Lu(T-t,\tilde X(t))dt.
$$
The first term in \eqref{e3.3} will be treated separately and we decompose the error as follows
\be
\label{e3.5}
u(T,x)-\E(\varphi(X_N)) =  u(T,x) -\E(u(T-\Dt, X_1))+\sum_{k=1}^{N-1} a_k +b_k +c_k.
\ee
Where
$$
\ba{l}
\ds a_k= \E \int_{t_k}^{t_{k+1}} \l(  A\tilde X(t)-A_\Dt X_k, Du(T-t,\tilde X(t))  \r)dt, \\
\\
\ds b_k= \E \int_{t_k}^{t_{k+1}} \l(   f(\tilde X(t))-S_\Dt f( X_k) , Du(T-t,\tilde X(t))  \r)dt,\\
\\
\ds c_k= \frac12 \E \int_{t_k}^{t_{k+1}}\Tr\l\{  \l[  \s(\tilde X(t))\s^*(\tilde X(t))
-\l( S_\Dt \s(X_k) \r)\l( S_\Dt \s(X_k) \r)^*   \r]  D^2 u(T-t,\tilde X(t)) \r\} dt.
\ea
$$
In the next steps, we estimate separately the different terms in \eqref{e3.5}.

\bigskip

{\bf Step 3:} Estimate of $u(T,x) -\E(u(T-\Dt, X_1))$.

\smallskip

By the Markov property
$$
u(T,x) = \E(\varphi(X(T,x)))= \E(u(T-\Dt,X(\Dt))).
$$
Therefore, by Lemma \ref{l1}, for any $\vep >0$,
$$
|u(T,x) -\E(u(T-\Dt, X_1))| \le c (T-\Dt)^{-1/2+\vep} \|\varphi\|_1\E \l( |X(\Dt)-X_1|_{-1/2+\vep}\r).
$$
Moreover
$$
\ba{ll}
X(\Dt)-X_1&= (S(\Dt) -S_\Dt)x + \int_0^\Dt S(t-s)f(X(s,x))ds  -\Dt S_\Dt f(x) \\
\\
&+\int_0^{\Dt} S(t-s)\s(X(s,x))dW(s)
-\sqrt{\Dt}S_\Dt \s(x) \chi_1.
\ea
$$
It is easy to prove that 
$$
|(-A)^{-1/2+\vep} \l(S(\Dt) -S_\Dt)\r|_{\L(H)}\le c \Dt^{1/2-\vep}.
$$
Since $(S(t))_{t\ge 0}$ is a contraction semigroup and $|(-A)^{-1/2+\vep}\cdot|\le c |\cdot|$, we have
by  \eqref{e2.3} and Lemma \ref{l5}
$$
\E \l|\int_0^\Dt S(t-s)f(X(s,x))ds\r|_{-1/2+\vep} \le  \Dt L_f \E(\sup_{s\in [0,\Dt]} |X(s,x)| +1)
\le c \Dt (|x|+1).
$$
Similarly
$$
|\Dt S_\Dt f(x)|_{-1/2+\vep} \le c \Dt (|x|+1).
$$
We then have
$$
\ba{l}
\ds \E\l(|\int_0^{\Dt} S(t-s) \s(X(s,x))dW(s) |_{-1/2+\vep}^2 \r) \\
\\
\ds = \E\l( \int_0^\Dt  | (-A)^{-1/2+\vep} S(t-s)\s(X(s,x))|_{\L_2(H)}^2 ds \r)\\
\\
\ds \le  \E\l( \int_0^\Dt  | (-A)^{-1/2+\vep}|_{\L_2(H)}^2 | S(t-s)|_{\L(H)}^2|\s(X(s,x))|_{\L(H)}^2 ds \r)\\
\ea
$$
and by  \eqref{e2.2},  \eqref{e2.4}, Lemma \ref{l5}
$$
\E\l(|\int_0^{\Dt} S(t-s) \s(X(s,x))dW(s) |_{-1/2+\vep}^2 \r)\le c \Dt (|x| +1). 
$$
Similarly
$$
\E\l( |\sqrt{\Dt}S_\Dt \s(x) \chi_1|^2\r) \le c \Dt (|x| +1).
$$
Gathering these estimate and using Cauchy-Schwartz inequality , we obtain
\be
\label{}
|u(T,x) -\E(u(T-\Dt, X_1))| \le c (T-\Dt)^{-1/2+\vep}  \Dt^{-1/2+\vep} \le c  \Dt^{1/2-\vep} 
\ee
where, as mentionned above, the constant is allowed to depend on $T$, $x$, $\varphi$, $f$, $\s$ $\dots$

\bigskip

{\bf Step 4:}  Estimate of $a_k$, $k\ge 1$.

\smallskip

We split $a_k$ as follows:
$$
a_k=a^1_k+a^2_k
$$
with
$$
\ba{l}
\ds a^1_k= \E \int_{t_k}^{t_{k+1}} \l( (A-A_\Dt)X_k, Du(T-t,\tilde X(t))  \r)dt,\\
\\
\ds a^2_k= \E \int_{t_k}^{t_{k+1}} \l( A(  \tilde X(t)- X_k), Du(T-t,\tilde X(t))  \r)dt.
\ea
$$
Note that $A_\Dt -A= \theta \Dt S_\Dt A^2$. By Lemma \ref{l1} below, we know that 
$Du(T-t,\tilde X(t))$ is in $D((-A)^\gamma)$ for $\gamma<1/2$ and it is easy to see that
$X_k$
belongs to $D((-A)^\delta)$ for $\delta<1/4$. It is impossible to compensate 
the presence of $A^2$ by such arguments. The idea is to recall \eqref{e2.7} and to 
observe that the irregularity of $X_k$ is contained in the stochastic integral. We thus 
further decompose $a^1_k$ in three terms according to \eqref{e2.7}. The first two terms 
are easy to treat. The third one involves the stochastic integral and is estimated thanks
to Malliavin calculus. We set
$$
\ba{l}
\ds a^{1,1}_k= -\theta \Dt\E \int_{t_k}^{t_{k+1}} \l( S_\Dt A^2 S_\Dt^k x, Du(T-t,\tilde X(t))  \r)dt,\\
\\
\ds a^{1,2}_k= -\theta \Dt\E \int_{t_k}^{t_{k+1}} \l( S_\Dt A^2\Dt \sum_{\ell=0}^{k-1} S_\Dt^{k-\ell} f(X_\ell) , Du(T-t,\tilde X(t))  \r)dt,\\
\\
\ds a^{1,3}_k=- \theta \Dt\E \int_{t_k}^{t_{k+1}} \l( S_\Dt A^2\sqrt\Dt \sum_{\ell=0}^{k-1} S_\Dt^{k-\ell}\s(X_{\ell})
\chi_{\ell +1 }, Du(T-t,\tilde X(t))  \r)dt,
\ea
$$ 
so that 
$$
a^1_k= a^{1,1}_k+a^{1,2}_k+a^{1,3}_k.
$$
By \eqref{e3}, \eqref{e1} and Lemma \ref{l1}, we have for $k=1,\dots,N-2$ and $\vep>0$
\be
\label{e3.6.1}
\ba{ll}
|a^{1,1}_k|&\ds \le c  \Dt \E\int_{t_k}^{t_{k+1}} | S_\Dt (-A)^{1/2+2\vep}|_{\L(H)}  | (-A)^{1-\vep} S^k_\Dt|_{\L(H)} 
|(-A)^{1/2-\vep}Du(T-t,\tilde X(t))| \; |x| dt \\
\\
&\ds \le c \Dt^{1/2-2\vep} t_k^{-1+\vep} \int_{t_k}^{t_{k+1}} (T-t)^{-(1/2-\vep)}dt.
\ea
\ee
The estimate of $a^{1,2}_k$ is similar. We have by \eqref{e2.3}, \eqref{e1}
$$
\ba{ll}
\ds \l|\Dt (-A)^{1-\vep} \sum_{\ell=0}^{k-1} S_\Dt^{k-\ell }f(X_{\ell})\r| &\ds \le 
L_f \Dt \sum_{\ell=0}^{k-1} \l| (-A)^{1-\vep}S_\Dt^{k-\ell}\r|_{\L(H)} \l(|X_{\ell}|+1\r)\\
\\
&\ds \le c  \Dt \sum_{\ell=0}^{k-1} t_{k-\ell}^{-1+\vep}  \l(|X_{\ell}|+1\r).
\ea
$$
Since
$$
\Dt \sum_{\ell=0}^{k-1} t_{k-\ell}^{-1+\vep} \le \vep^{-1} T^\vep,
$$
we deduce thanks to Lemma \ref{l1} and Lemma \ref{l2}
\be
\label{e3.6.3}
\ba{ll}
|a^{1,2}_k| &\le c \ds\Dt \int_{t_k}^{t_{k+1}} \l| S_\Dt (-A)^{1/2 + 2\vep} \r|_{\L(H)} (T-t)^{-1/2+\vep}dt\\
\\
&\le c \Dt^{1/2-2\vep} \int_{t_k}^{t_{k+1}} (T-t)^{-(1/2-\vep)} dt.
\ea
\ee
To treat $a^{1,3}_k$, we first rewrite it in terms of  a stochastic integral and then use 
Lemma \ref{l2.1}
$$
\ba{ll}
a^{1,3}_k&=\ds \theta \Dt\E \int_{t_k}^{t_{k+1}} \l(\int_0^{t_k} S_\Dt A^2 S_\Dt^{k-\ell_s} 
\s(X_{\ell_s})dW(s), Du(T-t,\tilde X(t))  \r)dt\\
\\
&=\ds \theta \Dt\E \int_{t_k}^{t_{k+1}} \int_0^{t_k} \Tr\l\{ \s^*(X_{\ell_s})
S_\Dt A^2 S_\Dt^{k-\ell_s}D^2u(T-t,\tilde X(t))D_s\tilde X(t) \r\} ds\; dt
\ea
$$
where $\ell_s=[s/\Dt]$ is the integer part of $s/\Dt$.
By the chain rule and \eqref{e3.0bis}, we have for $s\in [0,t_k]$, $h\in H$, $t\in [t_k,t_{k+1})$,
$$
D_s^h \tilde X (t) = D_s^h  X_k +  \int_{t_k}^{t} A_\Dt D_s^h X_k+  S_\Dt f'(X_k)\cdot D_s^h X_k ds +
  \int_{t_k}^{t}   S_\Dt \l(\s'(X_k)\cdot D_s^h X_k  \r) dW(s).
$$
Fro $\beta<1/4$, we have, by \eqref{e3},  \eqref{e2.2}, \eqref{e2.4}, \eqref{eq:HS}
$$
\ba{l}
\ds \E\l( \l| \int_{t_k}^{t}   S_\Dt \l(\s'(X_k)\cdot D_s^h X_k  \r) dW(s)\r|_\beta^2\r)\\
\\
\ds =\E\l(\int_{t_k}^t \l| (-A)^{\beta } S_\Dt  \l(\s'(X_k)\cdot D_s^h X_k  \r) \r|^2_{\L_2(H)} ds\r)\\
\\
\ds \le \E\l(\int_{t_k}^t \l|  (-A)^{-1/4-\vep}  \r|^2_{\L_2(H)} \l| (-A)^{\beta+1/4+\vep} S_\Dt \r|^2_{\L(H)}  \l|\s'(X_k)\cdot D_s^h X_k \r|^2_{\L(H)} ds\r)\\
\\
\le c \Dt^{1/2-2\beta-2\vep} \E\l( |D_s^h  X_k |^2\r).
\ea
$$
We then use \eqref{e3.0}, \eqref{e2.3}  to bound the other terms above and obtain thanks to Poincar\'e 
inequality
\be
\label{e3.6.4.1}
\E\l( |D_s^h \tilde X (t)|_\beta^2\r) \le c \E\l( |D_s^h  X_k |_\beta^2\r), \quad s\in [0,t_k],\quad t\in [t_k,t_{k+1}).
\ee
\bigskip
By Lemma \ref{l3}, we obtain for $\beta <1/4$
$$
\E\l( \l| (-A)^{\beta } D_s \tilde X (t)\r|_{\L(H)}^2\r) \le c t_{k-\ell_s}^{-2\beta }.
$$ 
We are now ready to conclude the estimate of $a_k^{1,3}$. We choose $\vep>0$ and write thanks to 
\eqref{e2.4}, \eqref{e3}, \eqref{e1}, Lemma \ref{l1bis} and  \eqref{e2.2}
$$
\ba{ll}
|a^{1,3}_k| \le & \ds \theta \Dt\E \int_{t_k}^{t_{k+1}}  \int_0^{t_k} 
\l|\s^*(X_{\ell_s})\r|_{\L(H)}
\l| S_\Dt A^{1/2+2\vep}\r|_{\L(H)}  \l| (-A)^{1-3\vep/2}  S_\Dt^{k-\ell_s}\r| _{\L(H)} \\
\\
&\times 
\ds \l| (-A)^{1/2-\vep/2} D^2u(T-t,\tilde X(t))(-A)^{1/2-\vep/2}\r|_{\L(H)} \Tr\l\{ (-A)^{-1/2-\vep/2}\r\}  \l| (-A)^{\vep} D_s\tilde X(t)\r|_{\L(H)}  ds\; dt\\
\\
&\ds \le c \Dt\E \int_{t_k}^{t_{k+1}}  \int_0^{t_k} \Dt^{-1/2-2\vep} t_{k-\ell_s}^{-1+3\vep/2} (T-t)^{-1+\vep}
t_{k-\ell_s}^{-\vep} ds\; dt.
\ea
$$
Since $\int_0^{t_k}t_{k-\ell_s}^{-1+\vep/2}  ds \le \frac2{\vep} T^{\vep/2}$, we deduce
\be
\label{e3.6.5}
|a^{1,3}_k| \le c \Dt^{1/2-2\vep}\int_{t_k}^{t_{k+1}}(T-t)^{-1+\vep}dt.
\ee
Gathering \eqref{e3.6.1}, \eqref{e3.6.3} and \eqref{e3.6.5}, we obtain for $k=1,\dots,N-1$
\be
\label{e3.7}
|a_k^1|\le c \Dt^{1/2-2\vep}\l( t_k^{-1+\vep}+1\r)\l(\int_{t_k}^{t_{k+1}}(T-t)^{-1+\vep}dt+1\r).
\ee

\bigskip

We now estimate  $a^2_k$.  Let us set
$$
\ba{l}
\ds a^{2,1}_k= \E\int_{t_k}^{t_{k+1}}(t-t_k) \l( AA_\Dt X_k,Du(T-t,\tilde X(t))\r)dt,\\
\\
\ds a^{2,2}_k= \E\int_{t_k}^{t_{k+1}}(t-t_k) \l( AS_\Dt f(X_k),Du(T-t,\tilde X(t))\r)dt,\\
\\
\ds a^{2,3}_k= \E\int_{t_k}^{t_{k+1}}\int_{t_k}^t \l( AS_\Dt \s(X_k)dW(s),Du(T-t,\tilde X(t))\r)dt,
\ea
$$
so that thanks to \eqref{e3.0bis}, we have
$a^2_k=a^{2,1}_k+a^{2,2}_k+a^{2,3}_k$. The first term $a^{2,1}_k$ is similar to $a^1_k$ above 
and is majorized in the same way
\be
\label{e3.9}
|a_k^{2,1}|\le c \Dt^{1/2-2\vep}\l( t_k^{-1+\vep}+1\r)\l(\int_{t_k}^{t_{k+1}}(T-t)^{-1+\vep}dt+1\r).
\ee
for $k=1,\dots,N-1$.
The second one is not difficult to treat, we have using similar arguments as above
\be
\label{e3.11}
\ba{ll}
\ds |a^{2,2}_k|&\ds \le c \Dt |(-A)^{1/2+\vep}S_\Dt|_{\L(H)} \E(|f(X_k)|)\int_{t_k}^{t_{k+1}} (T-t)^{-(1/2-\vep)}dt\\
\\
&\ds \le c \Dt^{1/2-\vep}\int_{t_k}^{t_{k+1}} (T-t)^{-(1/2-\vep)}dt
\ea
\ee
for $k=1,\dots,N-1$. The estimate of $a^{2,3}_k$ requires the use of Lemma \ref{l2.1}.  It implies
$$
a_k^{2,3}=\E \int_{t_k}^{t_{k+1}}\int_{t_k}^t \Tr\l\{  \s^*(X_k)S_\Dt A D^2u(T-t,\tilde X(t))D_s\tilde X(t) \r\}ds \;dt.
$$
Since, $X_k$ is $\F_{t_k}$ measurable, we have from \eqref{e3.0bis}
\be
\label{e3.11.1}
D_s \tilde X(t) = S_\Dt \s (X_k)\quad s\in (t_k,t_{k+1}],\quad t_k\le s\le t < t_{k+1}.
\ee
It follows, thanks to \eqref{e2.4}, \eqref{e3}, \eqref{e2.2} and Lemma \ref{l1bis},
$$
\ba{l}
a_k^{2,3}\ds =\E \int_{t_k}^{t_{k+1}}(t-t_k) \Tr\l\{  \s^*(X_k)S_\Dt A D^2u(T-t,\tilde X(t))S_\Dt \s (X_k) \r\} dt\\
\\
\ds \le c \Dt \E \int_{t_k}^{t_{k+1}} |\s(X_k)|_{\L(H)}|S_\Dt (-A)^{1/2+\vep/2} |_{\L(H)}
| (-A)^{1/2-\vep/2}D^2u(T-t,\tilde X(t))(-A)^{1/2-\vep/2} |_{\L(H)}  \\
\\
\hspace{6cm}\ds \times
\Tr( (-A)^{-1/2-\vep/2})| (-A)^{\vep}S_\Dt|_{\L(H)} |\s (X_k)|_{\L(H)} dt\\
\\
\ds \le c \Dt^{1/2-3\vep/2} \int_{t_k}^{t_{k+1}} (T-t)^{-1+\vep}dt
\ea
$$
for $k=1,\dots,N-1$ . Finally, we obtain
\be
\label{e3.13}
|a_k^{2}|\le c\Dt^{1/2-2\vep}(t_k^{-1+\vep}+1)\l(\int_{t_k}^{t_{k+1}} (T-t)^{-1+\vep}dt+1\r)
\ee
for $k=1,\dots,N-1$.
Together with \eqref{e3.7} this yields the estimate of $a_k$
$$
|a_k|\le c\Dt^{1/2-2\vep}(t_k^{-1+\vep}+1)\l(\int_{t_k}^{t_{k+1}} (T-t)^{-1+\vep}dt+1\r).
$$
It follows easily
\be
\label{e3.15}
\ds \sum_{k=1}^{N-1} |a_k|\le c \Dt^{1/2-2\vep}. 
\ee

\bigskip

{\bf Step 5:} Estimate of $b_k$.

\smallskip

This term seems easier to treat since we do not have the unbounded 
operator $A$. However, since it involves the nonlinear term, we need to use Ito formula 
\eqref{e3.1} to control $f(\tilde X(t))-f(X_k)$, this introduces many terms. For some of them we again use
Malliavin integration by parts.

First, we get rid of $S_\Dt$. We have thanks to \eqref{e4}, \eqref{e2.3}, Lemma \ref{l1} and 
Lemma \ref{l2}:
$$
\ba{ll}
b_k^1&\ds=\E\int_{t_k}^{t_{k+1}} \l(\l(I-S_\Dt\r) f(X_k),Du(T-t,\tilde X(t))\r)dt\\
&\ds\le c \E\int_{t_k}^{t_{k+1}} (1+|X_k|)|(-A)^{-1/2+\vep}(I-S_\Dt)|_{\L(H)}
|(-A)^{1/2-\vep}Du(T-t,\tilde X(t))|dt\\
&\ds\le c  \Dt^{1/2-\vep} \E\int_{t_k}^{t_{k+1}}(T-t)^{-1/2+\vep}dt
\ea
$$
for $k=0,\dots,N-1$.
We now estimate 
$$
\ba{ll}
b_k^2&\ds=b_k-b^1_k\\
&\ds=\E\int_{t_k}^{t_{k+1}} \l(f(\tilde X(t)-f(X_k),Du(T-t,\tilde X(t))\r)dt\\
&\ds =\E\int_{t_k}^{t_{k+1}} \sum_{i\in \N}( f_i(\tilde X(t)-f_i (X_k) )\partial_i u(T-t,\tilde X(t)) dt,
\ea
$$
where $f_i=(f,e_i)$ and $\partial_i = (D\cdot,e_i)$. We choose $(e_i)_{i\in\N}$ as the orthonormal basis of eigenvectors of $A$.
By \eqref{e3.1}, we have for $i\in\N$
$$
\ba{ll}
f_i(\tilde X(t) & \ds =f_i(X_k)+ \int_{t_k}^t \frac12\Tr\l\{  (S_\Dt\s(X_k))(S_\Dt\s(X_k))^* D^2f_i(\tilde X(s)) \r\}ds\\
\\
&\ds
+  \int_{t_k}^t  \l(A_\Dt X_k +S_\Dt f(X_k),Df_i(\tilde X(s))\r) ds +\int_{t_k}^t \l(Df_i(\tilde X(s)),\s(X_k)\r) dW(s).
\ea
$$
With obvious notations, this defines the decomposition  
$$
b_k^2=b_k^{2,1}+ b_k^{2,2}+b_k^{2,3}+b_k^{2,4}.
$$
To treat the first term, we rewrite it as follows{\footnote{Recall that we in fact work with Galerkin approximations so that all sums below are finite sums.}}:
$$
\ba{ll}
\ds b_k^{2,1}&=\ds \frac12\E\int_{t_k}^{t_{k+1}}\int_{t_k}^t \sum_{i\in \N}
\Tr\l\{  (S_\Dt\s(X_k))(S_\Dt\s(X_k))^* D^2f_i(\tilde X(s)) \r\}   \partial_i u(T-t,\tilde X(t)) ds \,dt\\
\\
& \ds = \frac12\E\int_{t_k}^{t_{k+1}}\int_{t_k}^t \Tr\l\{  (S_\Dt\s(X_k))(S_\Dt\s(X_k))^*  
{\mathcal A}(s,t)  \r\} ds \,dt
\ea
$$
where ${\mathcal A}(s,t)\in \L(H)$ is defined by
$$
\ba{ll}
({\mathcal A}(s,t)h,k)&\ds = \sum_{i\in \N}  D^2f_i(\tilde X(s)).(h,k)  \partial_i u(T-t,\tilde X(t))\\
& = \l(D^2f(\tilde X(s)).(h,k),
Du(T-t,\tilde X(t))\r),\; h,k\in H.
\ea
$$
Obviously
$$
\l|{\mathcal A}(s,t)\r|_{\L(H)}\le \l| D^2f(\tilde X(s))\r|_{\L^2(H\times H,H)} \l| Du(T-t,\tilde X(t)) \r|,
$$
where $\L^2(H\times H,H)$ denotes the space of  bilinear operators from $H\times H$ to $H$. By 
\eqref{e2.3} and Lemma \ref{l1}, we deduce:
$$
\l|{\mathcal A}(s,t)\r|_{\L(H)}\le c.
$$
Then, we write thanks to \eqref{e2.4}, \eqref{e3}, \eqref{e2.2},
$$
\ba{l}
\l|\Tr\l\{  (S_\Dt\s(X_k))(S_\Dt\s(X_k))^* {\mathcal A}(s,t) \r\}\r|\\
\le \Tr\l(\l(-A\r)^{-1/2-\vep}\r)
\l|(-A)^{1/2+\vep}S_\Dt\r|_{\L(H)}\l|\s(X_k)\r| ^2_{\L(H)}\l|{\mathcal A}(s,t)\r| _{\L(H)}\\
\le c \Dt^{-1/2-\vep}(1+|X_k|)^2.
\ea
$$
We deduce 
by  Lemma \ref{l2} 
\be
\label{}
b_k^{2,1}\le c \Dt^{3/2-\vep}.
\ee
The second term $b_k^{2,2}$ involves the same difficulty as $a^1_k$ above. We rewrite it using 
\eqref{e2.7}. This gives
$$
\ba{ll}
b_k^{2,2}&\ds = \E\int_{t_k}^{t_{k+1}}\int_{t_k}^t \sum_{i\in \N} \l(A_\Dt S_\Dt^k x + A_\Dt \Dt
\sum_{\ell=0}^{k-1}S_\Dt^{k-\ell} f(X_\ell), Df_i(\tilde X(s))\r) \partial_i u(T-t,\tilde X(t))ds dt\\
&\ds + \E\int_{t_k}^{t_{k+1}} \int_{t_k}^t \sum_{i\in \N} \l(A_\Dt \int_0^{t_k} 
S_\Dt^{k-\ell_\tau} \s(X_{\ell_\tau}) dW(\tau),  Df_i(\tilde X(s))\r) \partial_i u(T-t,\tilde X(t))ds dt
\ea
$$
where, as above, $\ell_\tau=[\tau/\Dt]$. The first term is bounded as follows, using Lemma \ref{l1}, \eqref{e2.3}, \eqref{e1}, \eqref{e3},
$$
\ba{l}
\ds \E\int_{t_k}^{t_{k+1}} \int_{t_k}^t \sum_{i\in \N} \l(A_\Dt S_\Dt^k x + A_\Dt \Dt
\sum_{\ell=0}^{k-1}S_\Dt^{k-\ell} f(X_\ell), Df_i(\tilde X(s))\r) \partial_i u(T-t,\tilde X(t))ds dt\\
\\
\ds = \E\int_{t_k}^{t_{k+1}} \int_{t_k}^t \l( Df(\tilde X(s))\cdot \l(A_\Dt S_\Dt^k x + A_\Dt \Dt
\sum_{\ell=0}^{k-1}S_\Dt^{k-\ell} f(X_\ell)\r), Du(T-t,\tilde X(t))\r)ds dt\\
\\
\ds \le c \E \int_{t_k}^{t_{k+1}} \int_{t_k}^t  \l| (-A)^\vep S_\Dt  \r|_{\L(H)} 
\bigg( \l| (-A)^{1-\vep}S_\Dt^k x  \r|  + 
\sum_{\ell=0}^{k-1} \l| (-A)^{1-\vep} S_\Dt^{k-\ell} \r|_{\L(H)} \l| f(X_\ell)\r| \bigg)ds\,dt\\
\\
\le c \Dt^{2-\vep}  ( t_k^{-1+\vep} +1).
\ea
$$
The second term of $b_k^{2,2}$ requires an integration by parts, we obtain
$$
\ba{l}
\ds \E\int_{t_k}^{t_{k+1}} \int_{t_k}^t \sum_{i\in \N} \l(A_\Dt \int_0^{t_k} 
S_\Dt^{k-\ell_\tau} \s(X_{\ell_\tau}) dW(\tau),  Df_i(\tilde X(s))\r) \partial_i u(T-t,\tilde X(t))ds dt\\
\\
= \ds  \E\int_{t_k}^{t_{k+1}} \int_{t_k}^t \sum_{i,j,m\in \N} \l(A_\Dt \int_0^{t_k} 
S_\Dt^{k-\ell_\tau} \s(X_{\ell_\tau})e_m,e_j\r) d\beta_m (\tau) \partial_j f_i(\tilde X(s))
 \partial_i u(T-t,\tilde X(t))ds dt\\
 \\
= 
 \ds  \E\int_{t_k}^{t_{k+1}} \int_{t_k}^t  \int_0^{t_k}  \sum_{i,j,m,n\in \N} \l(A_\Dt
S_\Dt^{k-\ell_\tau} \s(X_{\ell_\tau})e_m,e_j\r)  
\bigg[ \partial_{j,n} f_i(\tilde X(s))
\l(D_\tau^m\tilde X(s),e_n\r)  \partial_i u(T-t,\tilde X(t)) \\
\hfill +
\partial_j f_i(\tilde X(s))
 \partial_{i,n} u(T-t,\tilde X(t))\l(D_\tau^m\tilde X(t),e_n\r)\bigg]
d\tau ds dt\\
\ea
$$
$$
\ba{l}
= \ds  \E\int_{t_k}^{t_{k+1}} \int_{t_k}^t  \int_0^{t_k}  \sum_{i,m\in \N} D^2f_i(\tilde X(s)) 
\l(A_\Dt
S_\Dt^{k-\ell_\tau} \s(X_{\ell_\tau})e_m, D_\tau^m\tilde X(s)\r) \partial_i u(T-t,\tilde X(t))\\
\hfill + \l( B_i(s,t) A_\Dt
S_\Dt^{k-\ell_\tau} \s(X_{\ell_\tau})e_m,   D_\tau^m\tilde X(t)\r) d\tau ds dt\\
\\
= \ds  \E\int_{t_k}^{t_{k+1}} \int_{t_k}^t  \int_0^{t_k}  \sum_{i\in \N} \Tr \l\{ 
\l( D_\tau \tilde X(s)\r)^* D^2f_i(\tilde X(s)) 
A_\Dt
S_\Dt^{k-\ell_\tau} \s(X_{\ell_\tau})\r\}   \partial_i u(T-t,\tilde X(t))\hspace{1cm}\\
\hfill +   \Tr \l\{ \l( D_\tau \tilde X(t)\r)^*  B_i(s,t) A_\Dt
S_\Dt^{k-\ell_\tau} \s(X_{\ell_\tau} ) \r\}d\tau ds dt\\
\ea
$$
where, for $i\in\N$, $B_i(s,t)$ is defined by
$$
\l(B_i(s,t)g,h\r) \ds =(Df_i(\tilde X(s)),g)\sum_{n\in\N} \partial_{i,n} u(T-t,\tilde X(t))(h,e_n),\quad g,h\in H. 
$$
The first term above is estimate as $b_k^{2,1}$. For the second term, we write
$$
\ba{ll}
\sum_{i\in\N} \l(B_i(s,t)g,h\r)& = D^2u(T-t,\tilde X(t))\cdot ( Df(\tilde X(s))\cdot g,h) \\
\\
& = (D^2u(T-t,\tilde X(t))h, Df(\tilde X(s))\cdot g)         ,\quad g,h\in H. 
\ea
$$
Therefore
$$
\l|\sum_{i\in\N} B_i(s,t)\r|_{\L(H)} \le  \l|Df(\tilde X(s))\r|_{\L(H)}  \l|D^2u(T-t,\tilde X(t))\r|_{\L(H)} .
$$
We deduce by Lemma \ref{l3}, \eqref{e3.6.4.1}, \eqref{e2.3}, \eqref{e3}, \eqref{e2.2}, \eqref{e1},
\eqref{e2.2}, Lemma \ref{l2} and similar arguments as above
$$
\ba{l}
\ds \E\int_{t_k}^{t_{k+1}} \int_{t_k}^t \sum_{i\in \N} \l(A_\Dt \int_0^{t_k} 
S_\Dt^{k-\ell_\tau} \s(X_{\ell_\tau}) dW(\tau),  Df_i(\tilde X(s))\r) \partial_i u(T-t,\tilde X(t))ds dt\\
\\
\ds \le   c \Dt^{3/2-\vep}.
\ea
$$
Therefore
$$
b_k^{2,2} \le   c \Dt^{3/2-\vep}.
$$
It is also easy to see that 
$$
\ba{ll}
b_k^{2,3}&= \ds  \E\int_{t_k}^{t_{k+1}} \int_{t_k}^t \sum_{i\in \N}\l(S_\Dt f(X_k),Df_i(\tilde X(s))\r)
\partial_i u(T-t,\tilde X(t)) ds dt\\
&\ds = \E\int_{t_k}^{t_{k+1}} \int_{t_k}^t Du(T-t,\tilde X(t))\cdot \l( Df(\tilde X(s))\cdot S_\Dt f(X_k) \r) ds dt\\
\\
&\ds \le c \Dt^2.
\ea
$$
It remains to estimate $b_k^{2,4}$. We again integrate by parts the stochastic integral and obtain
by Lemma \ref{l2.1}:
$$
\ba{ll}
b_k^{2,4}&= \ds  \E\int_{t_k}^{t_{k+1}} \int_{t_k}^t \sum_{i\in \N}
\l( Df_i(\tilde X(s)),\s(X_k)dW(s)\r) \partial_i u(T-t,\tilde X(t)) dt\\
&\ds  = \E\int_{t_k}^{t_{k+1}} \int_{t_k}^t \Tr\l\{  \l(D_s\tilde X(t) \r)^* D^2 u(T-t,\tilde X(t)) 
Df(\tilde X(s))\s(X_k)\r\}  ds\,dt\\
& \ds = \E\int_{t_k}^{t_{k+1}} \int_{t_k}^t \Tr\l\{  \s^*(X_k) S_\Dt D^2 u(T-t,\tilde X(t)) 
Df(\tilde X(s))\s(X_k)\r\}  ds\,dt\\
&\ds \le c \Dt^{3/2-\vep},
\ea
$$
thanks to \eqref{e3.11.1}, \eqref{e2.2} and \eqref{e3}.

We conclude this step by gathering the previous estimates. This enables us to write
$$
\sum_{k=1}^{N-1} |b_k|\le c \Dt^{1/2-\vep}
$$

\bigskip

{\bf Step 6:} Estimate of $c_k$.

\smallskip

Using the symmetry of $Du$, we introduce the decomposition of $c_k$:

$$
\ba{ll}
c_k&\ds= \frac12 \E \int_{t_k}^{t_{k+1}}\Tr\l\{  \l[  \s(\tilde X(t))\s^*(\tilde X(t))
-\l( S_\Dt \s(X_k) \r)\l( S_\Dt \s(X_k) \r)^*   \r]  D^2 u(T-t,\tilde X(t)) \r\} dt\\
\\
&\ds= \frac12 \E \int_{t_k}^{t_{k+1}}\Tr\l\{  (I-S_\Dt)\s(\tilde X(t))\l((I-S_\Dt)\s(\tilde X(t))\r)^*
D^2 u(T-t,\tilde X(t)) \r\} dt\\
\\
&\ds+ \E \int_{t_k}^{t_{k+1}}\Tr\l\{   S_\Dt \s(\tilde X(t))\l((I-S_\Dt)\s(\tilde X(t))\r)^*
 D^2 u(T-t,\tilde X(t)) \r\} dt\\
 \\
 &\ds+ \frac12 \E \int_{t_k}^{t_{k+1}}\Tr\l\{  S_\Dt  \l(\s(\tilde X(t))-\s(X_k) \r)
\l( S_\Dt \s(\tilde X(t)) \r)^*     D^2 u(T-t,\tilde X(t)) \r\} dt\\
\\
&\ds+ \frac12 \E \int_{t_k}^{t_{k+1}}\Tr\l\{  S_\Dt \s(X_k) \l(S_\Dt \s(\tilde X(t))-\s(X_k) \r)^*     D^2 u(T-t,\tilde X(t)) \r\} dt\\
\\
&=c_k^1+c_k^2+c_k^3+c_k^4.
\ea
$$
The first two terms are easy to treat, we use similar arguments as in the previous steps and write
thanks to \eqref{eq:trace_prop2}, Lemma \ref{l1bis}, Lemma \ref{l2}, \eqref{e4}
$$
\ba{ll}
c_k^1& \le \ds c \E \int_{t_k}^{t_{k+1}} \Tr \l\{ (-A)^{-1/2+\vep  } (I-S_\Dt)
  \s(\tilde X(t) \s^*(\tilde X(t)) (I-S_\Dt)  (-A)^{-1/2+\vep  } \r\} (T-t)^{-1+2\vep} dt\\
\\
&\ds \le  c \E \int_{t_k}^{t_{k+1}} \Tr \l\{ (-A)^{-1/2+\vep  } (I-S_\Dt)
  (I-S_\Dt)  (-A)^{-1/2+\vep  } \r\} (T-t)^{-1+2\vep} dt\\
&\ds \le c \Dt^{1/2-3\vep} \int_{t_k}^{t_{k+1}} (T-t)^{-1+\vep} dt
\ea
$$
The second term is similar, we have
$$
\ba{ll}
c_k^2& \le \ds  \E \int_{t_k}^{t_{k+1}} \l|  (-A)^{-1/2+\vep  } (I-S_\Dt)\r|_{\L(H)} 
\l|  \s(\tilde X(t)\r|^2_{\L(H)}\l|  (-A)^{2\vep  } S_\Dt \r|_{\L(H)}  \Tr\l\{  (-A)^{-1/2-\vep  }   \r\} \\
&\hfill 
\l|  (-A)^{1/2-\vep  }D^2 u(T-t,\tilde X(t))(-A)^{1/2-\vep  }\r|_{\L(H)} dt\\
\\
&\ds \le c \Dt^{1/2-3\vep} \int_{t_k}^{t_{k+1}} (T-t)^{-1+\vep} dt
\ea
$$
The estimate of the  next term is much more complicated. It is based on similar arguments
as before but the computations are much longer. 

We use \eqref{e3.1} and obtain for $h,k\in H$:
$$
\ba{ll}
\l( \l( \s(\tilde X(t))-\s(X_k)\r)h,k\r) 
&=\ds \frac12\int_{t_k}^t\Tr \l\{  \l(S_\Dt\s(X_k)\r)\l(S_\Dt\s(X_k)\r)^* D^2\l( \s(\cdot)h,k\r) (\tilde X(s))\r\}dt\\
\\
& \ds + \frac12\int_{t_k}^t \l(  A_\Dt X_k +S_\Dt f(X_k),  D\l( \s(\cdot)h,k\r) (\tilde X(s)) \r)dt
\\
& \ds = ({\mathcal A}h,k)+({\mathcal B}h,k)+({\mathcal C}h,k).
\ea
$$
Thus we may write
$$
\ba{ll}
c^3_k & \ds = \frac12 \E \int_{t_k}^{t_{k+1}}\Tr\l\{  S_\Dt {\mathcal  A}
\l( S_\Dt \s(\tilde X(t)) \r)^*     D^2 u(T-t,\tilde X(t)) \r\} dt\\
\\
& +\ds  \frac12 \E \int_{t_k}^{t_{k+1}}\Tr\l\{  S_\Dt  {\mathcal B}
\l( S_\Dt \s(\tilde X(t)) \r)^*     D^2 u(T-t,\tilde X(t)) \r\} dt\\
\\
&\ds +  \frac12 \E \int_{t_k}^{t_{k+1}}\Tr\l\{  S_\Dt  {\mathcal C}
\l( S_\Dt \s(\tilde X(t)) \r)^*     D^2 u(T-t,\tilde X(t)) \r\} dt\\
\\
&= c_k^{3,1}+c_k^{3,2}+c_k^{3,3}.
\ea
$$
Note that
$$
({\mathcal A}h,k)= \frac12\int_{t_k}^t\sum_{\ell\in\N} \l(\l(\s''(\tilde X(s)).(S_\Dt \s(X_k)e_\ell,S_\Dt \s(X_k)e_\ell) \r)h,k\r) ds.
$$
By \eqref{e2.4bis}, for $u,v\in H$,
$$
\l(\l(\s''(\tilde X(s)).(u,v) \r)h,k\r)\le L_\s |u|_{-1/4}|v|_{-1/4}\, |h|\, |k| \le c |u|\, |v|\, |h|\, |k|.
$$
We deduce, thanks to \eqref{e2.2}, \eqref{e3},
$$
({\mathcal A}h,k) \le c \Dt^{1/2-\vep} (1+ |X_k|)^2 |h|\,|k|,
$$
and 
$$
|{\mathcal A}|_{\L(H)}\le c \Dt^{1/2-\vep}(1+ |X_k|)^2.
$$
Then, by Lemma \ref{l2}, Lemma \ref{l1bis}, \eqref{e3} and again \eqref{e2.2}
$$
c_k^{3,1} \le c \Dt^{1/2-3\vep} \int_{t_k}^{t_{k+1}} (T-t)^{-1/2+\vep}dt.
$$
The term $c_k^{3,2}$ involves the same difficulty as $a_k$ and $b_k^{2,2}$. We use 
\eqref{e2.7} to replace $X_k$ by a sum of three terms:
$$
\ba{ll}
({\mathcal B}h,k) &= \ds \frac12\int_{t_k}^t \bigg(  A_\Dt S_\Dt^k x + \Dt A_\Dt \sum_{\ell=0}^{k-1} S_\Dt^{k-\ell}f(X_\ell)
+ \int_0^{t_k} A_\Dt  S_\Dt^{k-\ell_\tau}\s(X_{\ell_\tau}) dW(\tau), \\
&\hfill  D\l( \s(\cdot)h,k\r) (\tilde X(s)) \bigg)ds\\
\\
&=\ds \frac12\int_{t_k}^t \bigg(\bigg[   \s'(\tilde X(s))\cdot\big( A_\Dt S_\Dt^k x + \Dt A_\Dt \sum_{\ell=0}^{k-1} S_\Dt^{k-\ell}f(X_\ell)\\
&\hfill
+ \int_0^{t_k} A_\Dt  S_\Dt^{k-\ell_\tau}\s(X_{\ell_\tau}) dW(\tau)  \big) \bigg]h,k\bigg)ds\\
\\
&= ({\mathcal B}_1h,k) + ({\mathcal B}_2h,k) +({\mathcal B}_3h,k) .
\ea
$$
We then write thanks to \eqref{e2.4}, \eqref{e1} and \eqref{e3}
$$
\ba{ll}
({\mathcal B}_1h,k) &= \ds \frac12\int_{t_k}^t \bigg(\bigg[   \s'(\tilde X(s))\cdot A_\Dt S_\Dt^k x  \bigg]h,k\bigg)ds\\
\\
& \ds \le c \int_{t_k}^t \l| \s'(\tilde X(s))\cdot A_\Dt S_\Dt^k x \r|_{\L(H)}|h|\, |k| ds\\
\\
&\le c \Dt |A_\Dt S_\Dt^k x |\, |h|\, |k|\\
\\
&\le c \Dt^{1-\vep} t_k^{1-\vep}\,|h|\, |k|  .
\ea
$$
Similarly
$$
({\mathcal B}_2h,k)\le c \Dt^{1-\vep} \,|h|\,|k|.
$$
It follows, thanks to Lemma \ref{l1bis}, \eqref{e3} and \eqref{e2.2}
$$
\ba{l}
\ds \frac12 \E \int_{t_k}^{t_{k+1}}\Tr\l\{  S_\Dt ( {\mathcal B}_1+{\mathcal B}_2)
\l( S_\Dt \s(\tilde X(t)) \r)^*     D^2 u(T-t,\tilde X(t)) \r\} dt\\
\\
\ds \le c \Dt^{1-3\vep} (t_k^{1-\vep}+1) \int_{t_k}^{t_{k+1}} (T-t)^{-1/2+\vep}dt.
\ea
$$

The estimate of the part of $c_k^{3,2}$ involving  ${\mathcal B}_3$ is very technical. As before, we get rid of the stochastic integral 
thanks to an integration by parts. This results in a supplementary trace term. In order to 
work with the double trace, we write everything in terms of the components of the operators 
and vectors. Given an operator $G$ on $H$, we set $G^{i,j}=(Ge_i,e_j)$. We thus 
write
$$
\ba{l}
\ds \E \int_{t_k}^{t_{k+1}}\Tr\l\{  S_\Dt  {\mathcal B}_3
\l( S_\Dt \s(\tilde X(t)) \r)^*     D^2 u(T-t,\tilde X(t)) \r\} dt\\
\\
\ds =  \E \int_{t_k}^{t_{k+1}} \sum_{i,j,m\in\N} {\mathcal B}_3^{i,j} \s^{m,j}(\tilde X(t))) \l(S_\Dt D^2u(T-t,\tilde X(t)) S_\Dt\r)^{m,i}dt\\
\ds =\sum_{i,j,m,n,r\in\N} \E \int_{t_k}^{t_{k+1}} \int_{t_k}^t \int_0^{t_k}  \partial_r\s^{i,j} (\tilde X (s)) \l(A_\Dt S_\Dt^{k-\ell_\tau}
\s(X_{\ell_\tau})e_n,e_r\r)d\beta_n(\tau) \s^{m,j}(\tilde X(t))) \\
\ds \hfill \l(S_\Dt D^2u(T-t,\tilde X(t)) S_\Dt\r)^{m,i} ds \, dt.
\ea
$$
It is important to recall here that in fact we work with finite dimensional approximations of the solutions so
that all the above sums are finite. We now use the Malliavin integration by parts and obtain
$$
\ba{l}
\ds \E \int_{t_k}^{t_{k+1}}\Tr\l\{  S_\Dt  {\mathcal B}_3
\l( S_\Dt \s(\tilde X(t)) \r)^*     D^2 u(T-t,\tilde X(t)) \r\} dt\\
\\
\ds =\sum_{i,j,m,n,r\in\N} \E \int_{t_k}^{t_{k+1}} \int_{t_k}^t \int_0^{t_k}  \sum_{p\in\N}
\partial_{r,p}\s^{i,j} (\tilde X (s))\l(   D_\tau^n \tilde X (s),e_p \r) \l(A_\Dt S_\Dt^{k-\ell_\tau}
\s(X_{\ell_\tau})e_n,e_r\r) \s^{m,j}(\tilde X(t))) \\
\ds \hfill \l(S_\Dt D^2u(T-t,\tilde X(t)) S_\Dt\r)^{m,i} \\
\\
\ds +  \sum_{p\in\N} \partial_r\s^{i,j} (\tilde X (s)) \l(A_\Dt S_\Dt^{k-\ell_\tau}
\s(X_{\ell_\tau})e_n,e_r\r)\partial_p \s^{m,j}(\tilde X(t))) \l( D_\tau^n \tilde X(t), e_p \r) 
 \l(S_\Dt D^2u(T-t,\tilde X(t)) S_\Dt\r)^{m,i} \\
\\
\ds +  \partial_r\s^{i,j} (\tilde X (s)) \l(A_\Dt S_\Dt^{k-\ell_\tau}
\s(X_{\ell_\tau})e_n,e_r\r) \s^{m,j}(\tilde X(t))) \l(S_\Dt \l(D^3 u(T-t,\tilde X(t))\cdot D_\tau^n \tilde X(s) \r)S_\Dt\r)^{m,i} 
d\tau \, ds\, dt\\
\\
= I+II+III.
\ea
$$
We then write
$$
\ba{l}
I=\ds  \sum_{i,j,m,n\in\N} \E \int_{t_k}^{t_{k+1}} \int_{t_k}^t \int_0^{t_k}
D^2\s^{i,j}(\tilde X(s))\cdot \l( D_\tau^n \tilde X (s),A_\Dt S_\Dt^{k-\ell_\tau}
\s(X_{\ell_\tau})e_n\r)\s^{m,j}(\tilde X(t))\\
\hfill \l(S_\Dt D^2u(T-t,\tilde X(t)) S_\Dt\r)^{m,i} d\tau \, ds\, dt\\
\\
= \ds  \sum_{j \in\N} \E \int_{t_k}^{t_{k+1}} \int_{t_k}^t \int_0^{t_k}
D^2u(T-t,\tilde X(t))\cdot\l(  \phi_1(s,\tau,k) e_j, S_\Dt \s(\tilde X(t))e_j \r) 
 d\tau \, ds\, dt\\
\\
=\ds  \E \int_{t_k}^{t_{k+1}} \int_{t_k}^t \int_0^{t_k} \Tr\l\{\s^*(\tilde X(t)) S_\Dt
D^2u(T-t,\tilde X(t))  \phi_1(s,\tau,k )\r\}
 d\tau \, ds\, dt
\ea
$$
where we have set
$$
\phi_1(s,\tau,k) h_1=\sum_{n\in \N} S_\Dt \l(D^2\s (\tilde X(s))\cdot \l( D_\tau^n \tilde X (s),A_\Dt S_\Dt^{k-\ell_\tau}
\s(X_{\ell_\tau})e_n\r)\r)h_1, \; h_1\in H.
$$
Let us define $\Sigma_{s,h_1,h_2}$ by
$$
\l( \Sigma_{s,h_1,h_2} u,v\r) = \l(S_\Dt \l(D^2\s (\tilde X(s))\cdot \l(u,v\r)\r) h_1,h_2\r), \; u,v\in H.
$$
Then by \eqref{e2.4bis}
$$
|\Sigma_{s,h_1,h_2}|_{\L(H)} \le c \, |h_1|\,|h_2|.
$$
We deduce by \eqref{e2.4}, \eqref{e2.2}, \eqref{e3}, \eqref{e1}, \eqref{e3.6.4.1} and Lemma \ref{l3}
$$
\ba{ll}
\l(\phi_1(s,\tau,k) h_1,h_2\r) & \ds = \Tr\l\{  \s^*(X_{\ell_\tau}) S_\Dt^{k-\ell_\tau} A_\Dt \Sigma_{s,h_1,h_2} 
D_\tau \tilde X (s) \r\}\\
\\
&\ds \le | \s^*(X_{\ell_\tau})|_{\L(H)} \Tr (-A)^{-1/2-\vep} 
| (-A)^{1/2+\vep}S_\Dt^{k-\ell_\tau} A_\Dt |_{\L(H)}  
|\Sigma_{s,h_1,h_2}|_{\L(H)} | D_\tau \tilde X (s)|_{\L(H)}\\
\\
&\le c \Dt^{-1/2-2\vep} t_k^{-1+\vep} |h_1|\,|h_2|(1+|X_k|)
\ea
$$
and by Lemma \ref{l2}, Lemma \ref{l1bis} and \eqref{e2.2}
$$
I\le c \Dt^{1/2-2\vep}t_k^{-1+\vep} \int_{t_k}^{t_{k+1}} (T-t)^{-1/2-\vep} dt.
$$
Similarly, we may write
$$
\ba{l}
II=\ds  \sum_{n\in \N} \E \int_{t_k}^{t_{k+1}} \int_{t_k}^t \int_0^{t_k} 
\Tr\bigg\{ \l[  \l(D\s(\tilde X (s))\cdot \l(A_\Dt S_\Dt^{k-\ell_\tau} \s(X_{\ell_\tau})e_n \r)\r)\r] \\
\\
\hfill \ds\l[\l( D\s(\tilde X(t))\cdot \l( D_\tau^n \tilde X(t) \r)\r)\r]^*
 S_\Dt D^2u(T-t,\tilde X(t))S_\Dt \bigg\} d\tau\, ds \, dt\\
 \\
 \ds \le c\Dt^{-2\vep} \E \int_{t_k}^{t_{k+1}} \int_{t_k}^t \int_0^{t_k} 
\l| \phi_2 (\tau,s,t,k)\r|_{\L(H)}
(T-t)^{-1/2+\vep} d\tau\, ds \, dt
 \ea
$$
with
$$
\phi_2(\tau,s,t,k)=\sum_{n\in\N} \l[  \l(D\s(\tilde X (s))\cdot \l(A_\Dt S_\Dt^{k-\ell_\tau} \s(X_{\ell_\tau})e_n \r)\r)\r] 
\l[ \l( D\s(\tilde X(t))\cdot \l( D_\tau^n \tilde X(t) \r)\r)\r]^*.
$$
We use similar arguments to estimate its norm. For $u,v\in H$, we write
$$
\ba{l}
(\phi_2(\tau,s,t,k)u,v)\\
\\
= \sum_{n\in\N}\l(
\l[\l( D\s(\tilde X(t))\cdot \l( D_\tau^n \tilde X(t) \r)\r)\r]^*u, 
 \l[  \l(D\s(\tilde X (s))\cdot \l(A_\Dt S_\Dt^{k-\ell_\tau} \s(X_{\ell_\tau})e_n \r)\r)\r]^* v\r)\\
 \\
 \ds=\Tr \l\{ \s^*(X_{\ell_\tau}) S_\Dt^{k-\ell_\tau} A_\Dt a_v^*  b_u   \r\} 
\ea
$$
with
$$
a_v h=  \l[ S_\Dt \l(D\s(\tilde X (s))\cdot  h \r)\r]^* v,
$$
$$
b_u h = \l[S_\Dt\l( D\s(\tilde X(t))\cdot \l( D_\tau^h \tilde X(t) \r)\r)\r]^*u.
$$
Since
$$
|a_v|_{\L(H)} \le c\, |v|, \quad |b_u|_{\L(H)} \le c\, |u|,
$$
we deduce 
$$
|\phi_2(\tau,s,t,k)|_{\L(H)} \le c \Tr \{ S_\Dt^{k-\ell_\tau} A_\Dt \} \le c \Dt^{-1/2-2\vep} t_{k-\ell_\tau}^{-1+\vep}. 
$$
and 
$$
II \le c \Dt^{1/2-4\vep}\int_{t_k}^{t_{k+1}} (T-t)^{-1/2-\vep}dt.
$$
Finally
$$
III= \frac12\int_{t_k}^{t_{k+1}} \int_{t_k}^{t}\int_0^{t_k}\sum_{n\in\N}
\Tr \{ \gamma_n S_\Dt \s(\tilde X(t))\} d\tau\, ds \, dt
$$
where for $u,v\in H$
$$
\ba{l}
\ds \sum_{n\in\N}(\gamma_nu,v)\\
\\
\ds= \sum_{n\in\N}D^3u(T-t,\tilde X(t))\l( D_\tau^n \tilde X(s), u, 
S_\Dt \l( D\s(\tilde X(s))\cdot\l( A_\Dt S_\Dt^{k-\ell_\tau} \s(X_{\ell_\tau}) e_n  \r) \r) v \r)\\
\\
=\Tr\l\{  \kappa(u,v) (-A)^{-1/2-\vep} (-A)^{2\vep} D_\tau^n \tilde X(s)  \r\},
\ea
$$
and for $h_1,h_2 \in H$
$$
( \kappa(u,v)h_1,h_2)= D^3u(T-t,\tilde X(t))\cdot \l((-A)^{1/2-\vep}h_1,u,  S_\Dt \l( D\s(\tilde X(s))\cdot\l( A_\Dt S_\Dt^{k-\ell_\tau}\s(X_{\ell_\tau})h_2  \r) \r) v\r).
$$
By Lemma \ref{l1ter} 
$$
|\kappa(u,v)|_{\L(H)}\le c(T-t)^{-1/2+\vep} t_{k-\ell_\tau}^{-1+3\vep} \Dt^{-3\vep} |u|\, |v|.
$$
Therefore, by \eqref{e2.2}, \eqref{e3.6.4.1} and Lemma \ref{l3},
$$
\ds \sum_{n\in\N}(\gamma_nu,v)\le c(T-t)^{-1/2+\vep} t_{k-\ell_\tau}^{-1+\vep} \Dt^{-3\vep} |u|\, |v|.
$$
It follows 
$$
|\gamma_n|_{\L(H)} \le c(T-t)^{-1/2+\vep} t_{k-\ell_\tau}^{-1+\vep} \Dt^{-3\vep}
$$
and by \eqref{e2.2}, \eqref{e3}
$$
III\le c   \Dt^{1/2-4\vep} \int_{t_k}^{t_{k+1}} (T-t)^{-1/2+\vep}dt.
$$
We can now conclude 
$$
 |c_k^{3,2} | \le c \Dt^{1/2-4\vep}(t_k^{-1+\vep}+1)(\int_{t_k}^{t_{k+1}} (T-t)^{-1/2+\vep}dt +1).
$$
Finally, it is easy to check 
$$
|{\mathcal C}|_{\L(H)} \le c \Dt (1+|X_k|)
$$
and 
$$
 |c_k^{3,3} | \le c \Dt^{1-2\vep}\int_{t_k}^{t_{k+1}} (T-t)^{-1/2+\vep}dt.
$$
We deduce
$$
 |c_k^{3} | \le c \Dt^{1/2-4\vep}(t_k^{-1+\vep}+1)\int_{t_k}^{t_{k+1}}( (T-t)^{-1/2-\vep} +1) dt,
$$
and, since $c_k^4$ is majorized in exactly the same way,
$$
 |c_k| \le c \Dt^{1/2-4\vep}(t_k^{-1+\vep}+1)\int_{t_k}^{t_{k+1}}( (T-t)^{-1+\vep} +1)dt.
$$
It follows
$$
\sum_{k=1}^{N-1}|c_k| \le c \Dt^{1/2-4\vep}.
$$

\bigskip

{\bf Step 7:} Conclusion.

\smallskip

It is now easy to gather all previous estimates in \eqref{e3.5} and 
deduce
$$
|u(T,x)-\E\l(\varphi(X_N)\r)|\le c \Dt^{1/2-4\vep}.
$$
Recall that all the above computations have been done on the Galerkin approximations of 
$X$ and $X_k)$. The constant $c$ above does not depend on $m$ so that we can 
easily let $m\to \infty $ in this estimate and obtain the result.

\section{Auxiliary Lemmas}

In this section, we state and prove technical Lemmas used in the preceeding section. Again, the various
estimates used here could be difficult to justify rigorously on the infinite dimensional equation and we in
fact work with Galerkin approximations. Taking the limit $m\to \infty$ at the end of the proofs gives the
results rigorously.

The first two Lemmas are very classical and we state them without proof.
\begin{Lemma}
\label{l2}
For any $l\in \N$, there exists a constant $c_l$ such that
$$
\max_{k=0,\dots, N} \E(|X_k|^l) \le c_l (|x|^l+1).
$$
\end{Lemma}


\begin{Lemma}
\label{l5}
For any $l\in \N$, there exists a constant $\tilde c_l$ such that
$$
\sup_t \E(|X(t,x)|) \le \tilde c_l (|x|^l+1).
$$
\end{Lemma}

\begin{Lemma}
\label{l3}
For any $\beta \in [0,1/4)$, there exists a constant $c$ such 
for  $k=1,\dots,N$ , $s\in [0,t_k]$, we have
$$
t_{k-\ell_s}^{2\beta } \E\l(|D_s^h X_k|_\beta^2\r) \le c |h|^2, \; h\in H.
$$
\end{Lemma}
{\bf Proof:} By \eqref{e2.7} and the chain rule, we obtain the following formula
for the Malliavin derivative of $X_k$:
$$
\ba{ll}
D_s^h X_k&= \ds
 S_\Dt^{k-\ell_s} \s(X_{\ell_s})h +\Dt \sum_{\ell=\ell_s+1}^{k-1} S_\Dt^{k-\ell}  f'(X_{\ell}) 
\cdot D_s^hX_{\ell}\\
\\
&+ 
\ds\sqrt\Dt\sum_{\ell=\ell_s+1}^{k-1} S_\Dt^{k-\ell} (\s'(X_{\ell})\cdot D_s^hX_{\ell} )\chi_{\ell+1}
\ea
$$
for $s\in [0,t_k]$ and $h\in H$.

By \eqref{e1}, \eqref{e2.2},\eqref{e2.3}, \eqref{e2.4}, we deduce for $\vep>0$
$$
\ba{lcl}
\ds\E\l(|D_s^h X_k|_\beta^2\r)& \le &\ds c\bigg( t_{k-\ell_s}^{-2\beta } |h|^2 +  
\l(\Dt  \sum_{\ell=\ell_s+1}^{k-1} \l|(-A)^\beta S_\Dt^{k-\ell}\r|_{\L(H)}  \l|f'(X_{\ell})\r|_{\L(H)} \l|D_s^hX_{\ell}\r| 
\r)^2\\
\\
&&\ds+\Dt  \sum_{\ell=\ell_s+1}^{k-1} \l|(-A)^\beta S_\Dt^{k-\ell} (\s'(X_{\ell})\cdot D_s^hX_{\ell} )\r|_{\L_2(H)} ^2\bigg) \\
\\
&\le & \ds c\bigg( t_{k-\ell_s}^{-2\beta } |h|^2 +  
L_F^2 \l(\Dt  \sum_{\ell=\ell_s+1}^{k-1} t_{k-\ell}^{-\beta}  \l|D_s^hX_{\ell}\r| 
\r)^2\\
\\
&&\ds+L_\s \Dt  \sum_{\ell=\ell_s+1}^{k-1}t_{k-\ell}^{-1/2-\vep-2\beta}  \l|D_s^hX_{\ell} )\r|_{\L_2(H)} ^2\bigg).
\ea
$$
It is now easy to use a discrete Gronwall Lemma and prove
$$
\ds\max_{l=\ell_s+1,\dots,k} t_{\ell-\ell_s}^{2\beta } \E\l(|D_s^h X_\ell|^2\r) \le c |h|^2
$$ 

\hfill $\square$

\begin{Lemma}
\label{l1} Let $\varphi \in C^1_b(H,\R)$. For any $\beta  <1/2$, there exists a constant $c_\beta $ such that for 
$t>0$, $x\in H$
$$
|Du(t,x)|_\beta \le c_\beta  t^{-\beta } \|\varphi \|_{1},
$$
where $u$ is defined in \eqref{e3.2}.
\end{Lemma}
{\bf Proof:}  Differentiating \eqref{e3.2}, we obtain for $h\in H$:
$$
Du(t,x)\cdot h = \E\l( D\varphi(X(t,x))\cdot \eta^{h,x}(t) \r)
$$
where $\eta^{h,x}(t)$ is the solution of 
$$
\l\{
\ba{l}
d\eta^{h,x} = \l(   A\eta^{h,x}+f'(X(t,x))\cdot \eta^{h,x} \r)dt + \s'(X(t,x))\cdot \eta^{h,x} dW,\\
\\
\eta^{h,x}(0)=h.
\ea
\r.
$$
We rewrite this equation in the integral form 
$$
\eta^{h,x}(t)=S(t)h+\int_0^t S(t-s) f'(X(s,x))\cdot \eta^{h,x} (s) ds + \int_0^t S(t-s)\s'(X(s,x))\cdot \eta^{h,x}(s) dW(s),\; t\ge 0.
$$
By \eqref{e2.4}, \eqref{e2.2}, \eqref{e1}, we have for $y,k\in H$ and $\alpha >1/2$:
$$
\l|S(t) \s'(y)\cdot k\r|_{\L_2(H)} \le L_\s  \l| (-A)^{-\alpha/2}\r|_{\L_2(H)} \l|(-A)^{\alpha/2}S(t)\r|_{\L(H)}  |k|
\le c t^{-\alpha/2} |k|.
$$
Using \eqref{e2.3} and then Cauchy-Schwarz inequality, we obtain
$$
\ba{ll}
\ds \E\l(\l|\eta^{h,x}(t)\r|^2\r)& \ds \le c \; t^{-2\beta }  |h|_{-\beta }^2
+L_f^2 \E\l(\l(\int_0^t  |\eta^{h,x} (s)| ds\r)^2\r) 
+ \E\int_0^t (t-s)^{-\alpha}\l|\eta^{h,x}(s)\r|^2 ds\\
\\
& \ds \le c \; t^{-2\beta }  |h|_{-\beta }^2
+c \int_0^t  \E(|\eta^{h,x} (s)|^2 ds  
+ \E\int_0^t (t-s)^{-\alpha}\l|\eta^{h,x}(s)\r|^2 ds.
\ea
$$
It is classical that this implies 
\be
\label{a1}
\sup_{t\in [0,T]} t^{2\beta} \E\l(\l|\eta^{h,x}(t)\r|^2\r) \le |h|_{-\beta }^2.
\ee
We deduce 
$$
|Du(t,x)\cdot h |\le c \| \varphi\|_{1}  t^{-\beta}|h|_{-\beta }.
$$
Taking the supremum over $h$ yields the result.

\hfill $\square$

\begin{Lemma}
\label{l1bis} Let $\varphi \in C^2_b(H,\R)$. For any $\beta,\gamma <1/2$, there exists a constant $c_{\beta,\gamma} $ such that for 
$t>0$, $x\in H$
$$
|(-A)^{\beta}D^2u(t,x)(-A)^{\gamma}|_{\L(H)} \le c_{\beta,\gamma}  t^{-(\beta+\gamma) } \|\varphi \|_{2},
$$
where $u$ is defined in \eqref{e3.2}.
\end{Lemma}
{\bf Proof:} We use the same notations as in the proof of Lemma \ref{l1}. We differentiate
a second time  \eqref{e3.2} and obtain for $h,k\in H$:
\be
\label{a1bis}
D^2u(t,x)\cdot (h,k) = \E\l( D^2\varphi(X(t,x))\cdot (\eta^{h,x}(t),\eta^{k,x}(t))+ 
D\varphi(X(t,x))\cdot \zeta^{h,k,x}(t) \r)
\ee
where $\zeta^{h,k,x}(t)$ is the solution of 
$$
\l\{
\ba{ll}
d\zeta^{h,k,x} &= \l(   A\zeta^{h,k,x}+f''(X(t,x))\cdot (\eta^{h,x}(t),\eta^{k,x}(t))
+f'(X(t,x))\cdot \zeta^{h,k,x}(t)\r)dt \\
\\
&+ \l(\s''(X(t,x))\cdot (\eta^{h,x}(t),\eta^{k,x}(t)) 
+\s'(X(t,x))\cdot \zeta^{h,k,x}(t) \r)dW,\\
\\
\zeta^{h,k,x}(0)&=0.
\ea
\r.
$$
We rewrite this equation in the integral form 
$$
\ba{ll}
\zeta^{h,k,x}(t)&\ds=\int_0^t S(t-s)\l(f''(X(s,x))\cdot (\eta^{h,x}(s),\eta^{k,x}(s))+
 f'(X(s,x))\cdot \zeta^{h,k,x} (s)\r) ds\\
\\
&\ds + \int_0^t S(t-s) \l(\s''(X(s,x))\cdot (\eta^{h,x}(s),\eta^{k,x}(s)) 
+\s'(X(s,x))\cdot \zeta^{h,k,x}(s) \r) dW(s),\; t\ge 0.
\ea
$$
Using similar argument as above and \eqref{e2.4bis}, we prove
\be
\label{a2}
\ba{ll}
\E\l(|\zeta^{h,k,x}(t)|^2\r) &\ds\le c \E \l(\int_0^t |\eta^{h,x}(s)|\, |\eta^{k,x}(s)|+ |\zeta^{h,k,x} (s)|  ds\r)^2\\
\\
&\ds + c\E \int_0^t  (t-s)^{-\alpha}\l(|\eta^{h,x}(s)|_{-1/4}^2 |\eta^{k,x}(s)|_{-1/4}^2 
+|\zeta^{h,k,x}(s)|^2 \r) ds,\; t\ge 0.
\ea
\ee
Proceeding as in Lemma \ref{l1}, we have thanks to Burkholder inequality and then to 
Minkowsky inequality
$$
\ba{ll}
\ds \E\l(\l|\eta^{h,x}(t)\r|^4\r)& \ds \le c \; t^{-4\beta }  |h|_{-\beta }^4
+c \E\l(\l(\int_0^t  |\eta^{h,x} (s)| ds\r)^4\r) 
+ \E\l( \l(\int_0^t (t-s)^{-\alpha}\l|\eta^{h,x}(s)\r|^2 ds\r)^2\r)\\
\\
& \ds \le c \; t^{-4\beta }  |h|_{-\beta }^4
+c \l(\int_0^t  \l(\E(|\eta^{h,x} (s)|^4)\r)^{1/4} ds\r)^{4}  \\
\\
&\ds+ c\l( \int_0^t (t-s)^{-\alpha}\l(\E(|\eta^{h,x} (s)|^4)\r)^{1/2} ds\r)^{2}.
\ea
$$
Taking the square root of this inequality and using a generalized Gronwall Lemma, we deduce 
\be
\label{a3}
\sup_{t\in [0,T]} t^{4\beta} \E\l(\l|\eta^{h,x}(t)\r|^4\r) \le c |h|_{-\beta }^4.
\ee
Similarly, we have
$$
\ba{ll}
\E(|\eta^{h,x}(t)|_{-1/4}^4)&\ds \le c t^{1-4\beta}  |h|_{-\beta }^4
+ c\E \l(\int_0^t |\eta^{h,x}(s)| ds\r)^4 + c\E \l(\int_0^t (t-s)^{-\alpha}|\eta^{h,x}(s)|^2 ds\r)^2\\
&\ds \le c t^{1-4\beta}  |h|_{-\beta }^4
+  c\l(\int_0^t \E\l(|\eta^{h,x}(s)|^4\r)^{1/4} ds\r)^4\\
& \ds+c  \l(\int_0^t (t-s)^{-\alpha}\E\l(|\eta^{h,x}(s)|^4\r)^{1/2} ds\r)^2
\ea
$$
Therefore, by \eqref{a3},
\be
\label{a4}
\ba{ll}
\E(|\eta^{h,x}(t)|_{-1/4}^4)&\ds \le  c t^{1-4\beta}  |h|_{-\beta }^4
\ea
\ee

Plugging these inequalities and similar ones for $\eta^{k,x}$ in \eqref{a2} yields
$$
\sup_{t\in [0,T]}  \E\l(\l|\zeta^{h,x}(t)\r|^2\r) \le c |h|_{-\beta }^2|h|_{-\gamma }^2.
$$
The result follows easily using \eqref{a3} and this inequality in \eqref{a1bis}. 

\hfill $\square$

The following Lemma is proved thanks to similar arguments.

\begin{Lemma}
\label{l1ter}
Let $\varphi \in C^3_b(H,\R)$. For any $\beta <1/2$, there exists a constant $c_{\beta} $ such that for 
$t>0$, $x\in H$, $h_1\in D((-A)^{\beta})$, $h_2\in H$, $h_3\in H$
$$
D^3u(t,x)\cdot ((-A)^{\beta}h_1, h_2, h_3) \le c_{\beta}  t^{-\beta } \|\varphi \|_{3} |h_1|\, |h_2|\, |h_3|,
$$
where $u$ is defined in \eqref{e3.2}.

\end{Lemma}

\end{document}